\documentclass[opre,nonblindrev]{informs3} 

\OneAndAHalfSpacedXII 

\usepackage[ruled,noend]{algorithm2e}
\usepackage{subfig}
\usepackage{scalerel}
\usepackage{mathrsfs}
\usepackage{mathtools}
\usepackage{bbm}

\usepackage{natbib}
 \bibpunct[, ]{(}{)}{,}{a}{}{,}%
 \def\bibfont{\small}%
\TheoremsNumberedThrough     

\EquationsNumberedThrough    

\allowdisplaybreaks

\begin{document}


\RUNAUTHOR{Yang, Nadarajah, Secomandi}

\RUNTITLE{Pathwise Optimization for Merchant Energy Production}

\TITLE{Pathwise Optimization for Merchant Energy Production}

\ARTICLEAUTHORS{%
\AUTHOR{Bo Yang}
\AFF{Tepper School of Business, Carnegie Mellon University, Pittsburgh, PA, USA, \EMAIL{boy1@andrew.cmu.edu}}
\AUTHOR{Selvaprabu Nadarajah}
\AFF{College of Business, University of Illinois at Chicago, Chicago, IL, USA, \EMAIL{selvan@uic.edu}}
\AUTHOR{Nicola Secomandi}
\AFF{Tepper School of Business, Carnegie Mellon University, Pittsburgh, PA, USA, \EMAIL{ns7@andrew.cmu.edu}}
} 

\ABSTRACT{%
{
We study merchant energy production modeled as a compound switching and timing option.
The resulting Markov decision process is intractable.
State-of-the-art approximate dynamic programming methods applied to realistic instances of this model yield policies with large optimality gaps that are attributed to a weak upper (dual) bound on the optimal policy value. We extend pathwise optimization from stopping models to merchant energy production to investigate this issue. We apply principal component analysis and block coordinate descent in novel ways to respectively precondition and solve the ensuing ill conditioned and large scale linear program, which even a cutting-edge commercial solver is unable to handle directly. Compared to standard methods, our approach leads to substantially tighter dual bounds and smaller optimality gaps at the expense of considerably larger computational effort. Specifically, we provide numerical evidence for the near optimality of the operating policies based on least squares Monte Carlo and compute slightly better ones using our approach on a set of existing benchmark ethanol production instances. These findings suggest that both these policies are effective for the class of models we investigate. Our research has potential relevance for other commodity merchant operations settings.
}
}%

\HISTORY{December 2019}

\maketitle

%

\section{Introduction}
We study the merchant management of energy production assets, such as power and natural-gas-processing plants, oil and bio refineries, and ethanol manufacturing facilities \citep{ tseng2002short, tseng2007framework, devalkar2011integrated,  thompson2013optimal, dong2014value, kn:BoyabatliEtAl2017,  nadarajah2018least}.
Modeling as a portfolio of real options the ability of the managers of these assets to dynamically adapt their operating policies to changing market conditions provides a convenient approach to maximize their market values \citep{kn:DixitPindyck1994,trigeorgis1996real,kn:SmithMcCardle1998,smith1999options, kn:EydelandWolyniec2003,kn:Geman2005,smith2005alternative,guthrie2009real,secomandi2014real,kn:SecomandiSeppi2016}.


Managing an ethanol factory in wholesale markets \citep[Chapter 17]{guthrie2009real} exemplifies the main ideas underlying merchant energy production.
Operating the plant is desirable when the spread between the output and input wholesale prices net of the conversion cost is attractive.
Temporary or prolonged periods of unappealing spreads can be dealt with by suspending production or mothballing the plant.
In the latter case, reactivating or abandoning operations is advisable when the spread improves or worsens.
This managerial flexibility can be modeled as a compound switching and timing option on the uncertain evolution of the prices of ethanol and of corn and natural gas (the raw materials).

As is typical in merchant energy operations \citep{secomandi2014real, kn:SecomandiSeppi2016}, real option models of energy production give rise to intractable Markov decision processes (MDPs).
In each stage the MDP state contains the status of both the plant and the market.
The choices of the merchant producer determine the evolution of the former component.
Given stochastic processes govern the dynamics of the latter one and, based on price taking and small plant size assumptions, are independent of these decisions.
Intertemporal linkages between operational conditions and high dimensional market information (input and output current futures curves) lead to (some of) the well-known ``curses of dimensionality" (\citealt{kn:Powell2011}, \S 1.2).
Approximate dynamic programming (ADP) methods are thus used to obtain (operating) policies and optimality gaps.

Combining least squares Monte Carlo (LSM) and information relaxation and duality techniques (\citealp{carriere1996valuation, longstaff2001valuing, smith2005alternative, cortazar2008valuation, brown2010information, kn:SecomandiSeppi2016, nadarajah2017comparison, kn:Secomandi2017}) is a state-of-the-art ADP approach for intractable merchant operations MDPs.
LSM computes a value function approximation (VFA) that defines both a policy and penalties on hindsight information, which in turn respectively lead to estimates of lower and upper (dual) bounds on the optimal policy value.
\citet{nadarajah2018least} observe that this procedure yields a sizable (about 12\%) average optimality gap on realistic merchant ethanol production instances.
They attribute this finding primarily to the looseness of their estimated dual bounds.



Pathwise optimization (PO, \citealp{desai2012pathwise, chandramouli2019convex}) is an ADP approach that first solves a linear program formulated on a set of Monte Carlo sample paths of the market uncertainty to find the best dual penalties constructed using VFAs that, as in LSM, are linear combinations of basis functions.
It then computes a corresponding policy and optimality gap.
We extend this methodology from stopping models \citep{desai2012pathwise, chandramouli2019convex} to merchant energy production to assess the quality of the LSM effectiveness.

Our PO linear program (PLP) specified using the basis functions used by \citet{nadarajah2018least}, which are common in the literature (see, e.g., \citealp{nadarajah2017comparison} and references therein), is hard to solve even using a state-of-the-art commercial solver (Gurobi). This difficulty occurs even for instances that are significantly smaller than the ones considered by these authors, because PLP is both ill conditioned and large scale.
These drawbacks respectively arise from the considerable parallelism of the chosen basis functions and the need to represent the dual value function for each stage and state as a decision variable to ensure that the PLP size grows linearly with the number of stages.
We address these issues by developing (i) an exact preconditioning procedure based on principal component analysis (PCA) and (ii) a block coordinate descent (BCD) optimization method.
Our PCA procedure exploits the block diagonal structure of PLP for efficiency and exactly reformulates it by orthogonalizing its columns within each block.
Our BCD algorithm solves this linear program by iteratively and cyclically optimizing the values of groups of decision variables, while fixing the ones of the remaining variables.
It thus requires less memory than employing a monolithic approach, which is impractical due to its excessive memory requirement. We establish that an idealized version of our BCD technique that satisfies common assumptions converges to an optimal PLP solution.

We apply to the instances of \citet{nadarajah2018least} LSM and our PO approach, both specified with VFAs expressed only for states that admit more than one action.
The PO-based lower and dual bound estimates respectively outperform slightly and substantially the LSM-based ones, by 1\% and 5\% on average.
The respective estimated optimality gaps of the PO- and LSM-based policies, both obtained using the PO-based dual bound, are on average 7\% and 8\%, which compare favorably with the corresponding 12\% LSM-based figure obtained by \cite{nadarajah2018least}.
These results provide evidence for the effectiveness of both these policies for the type of models that we study.
PO is considerably more computationally onerous than LSM, on average taking eleven hours instead of seven minutes per instance. However, our algorithms allow us to find high quality solutions to the resulting linear programs, which cannot be otherwise optimized.


Our research is potentially relevant for other commodity merchant operations contexts and related real option models \citep{secomandi2014real, kn:SecomandiSeppi2016}, including oil and natural gas extraction fields, liquefied natural gas facilities, copper mines, and renewable energy plant \citep{brennan1985evaluating, kn:SmithMcCardle1998, smith1999options, smith2005alternative, cortazar2008valuation,romo2009optimizing, enders2010interaction, kn:LaiEtAl2011, arvesen2013linepack, denault2013simulation,kn:HinzYee2018,kn:ZhouEtAl2019}.
The use of PCA as a preconditioning technique may have applicability in mathematical programming beyond our specific application. Because we use PLP to represent a piecewise linear and convex objective function, opportunely modified variants of our BCD algorithm and the convergence analysis of its idealized version may apply to other models that optimize analogous functions.

We discuss the novelty of our contributions in \S\ref{Sec:LitRev}.
Section \ref{Sec:Model} presents the merchant energy production MDP that we study.
We introduce PO, formulate PLP, and describe how to obtain both a policy and optimality gaps from its solution in \S\ref{Section 3}.
Section \ref{Section 4} explains how we use PCA and BCD to solve PLP.
We report the results of our numerical investigation in \S\ref{Section 5}.
Section \ref{Sec:Conc} concludes. An appendix includes the proofs of formal results.

\section{Literature Review}
\label{Sec:LitRev}

Our VFA representation is more parsimonious than the one of \citet{nadarajah2018least} because it excludes states with fixed actions.
\citet{kn:Trivella2019} study a version of their model that discourages abandonment and adopt LSM and duality methods.
They observe that on average the LSM-based bounds are tighter in their application than they are in ours.

PLP differs from the PO linear programs of \citet{desai2012pathwise} and \citet{chandramouli2019convex} for single and multiple stopping models, respectively, because we represent the dual value functions as decision variables, whereas they respectively enumerate the individual or cumulative penalized payoff of every possible stopping choice or sequence of such alternatives.
\cite{desai2012pathwise} apply a commercial linear programming solver to readily optimize their model. \cite{chandramouli2019convex} approximately solves his model as a sequence of single stopping models formulated and solved as in \cite{desai2012pathwise}. In contrast, our PLP solution approach relies on PCA and BCD methods.
Our findings on the performance of the LSM- and PO-based policies and dual bounds are largely consistent with the ones of these authors. 



Preconditioning is a common technique to facilitate solving mathematical programs, in particular with linear constraints (see, e.g., \citealp{renegar1995incorporating, renegar1995linear, cheung2001new, epelman2002new, belloni2009geometric, amelunxen2012coordinate, pena2014some}, and references therein).
Our use of PCA for preconditioning, rather than the typical dimensionality reduction purpose (see, e.g., \citealp{kn:Jolliffe2002}), appears to be unique. The literature on BCD algorithms for optimization is extensive (see, e.g., \citealp{sargent1973convergence, grippo2000convergence, PTseng2001, nesterov2012efficiency, richtarik2014iteration, bertsekas2015convex} Chapter 6, and references therein). Our use of BCD in a PO setting is novel. The idealized version of our BCD algorithm and its theoretical analysis rely on common assumptions (see, e.g., \citealt{bertsekas2015convex}, \S 6.5)
\section{Model}
\label{Sec:Model}

We formulate an MDP for the merchant management of energy production modeled as a compound switching and timing option.
This material is largely from \citet[\S7.1]{nadarajah2018least}, which is based on \citet[Chapter 17]{guthrie2009real}.

For concreteness, consider a plant that converts corn and natural gas into ethanol (with minor changes, the MDP that we formulate applies to other energy production settings).
This facility is managed on a merchant basis in wholesale markets for these commodities.
Operating choices are made periodically over a finite horizon in the face of uncertain input and output prices.
Specifically, if the plant is operational, the merchant can produce at full capacity, suspend production, start the mothballing (temporary abandonment) process, or abandon operations.
Upon completion of the mothballing activity, the facility can be kept mothballed, abandoned, or reactivated by initiating the reactivation process, at the end of which the plant becomes operational.
Abandonment is compulsory at the end of the horizon, unless it has occurred earlier.
No interruption of the mothballing and reactivation processes is possible.

The finite horizon consists of $I$ decision dates.
The stage set $\mathcal{I}:= \{0,1,\ldots,I-1\}$ includes their indices.
The evolution of the input and output market prices is Markovian.
Each state features both the plant operating mode and the forward curves (defined below) of the inputs and the output.

The operational and abandoned operating modes are respectively $\mathsf{O}$ and $\mathsf{A}$.
The respective durations of the mothballing and reactivation processes are $N^\mathsf{M}$ and $N^\mathsf{R}$ stages.
The set $\mathcal{M}$ contains the partially and fully mothballed operating modes $\mathsf{M}_{1}$ through $\mathsf{M}_{N^{\mathsf{M}}-1}$ and $\mathsf{M}_{N^{\mathsf{M}}}$.
We include in set $\mathcal{R}$ the partially reactivated operating modes ${\mathsf{R}}_{1}$ through ${\mathsf{R}}_{N^{\mathsf{R}} - 1}$.
We define the stage $i$ operating mode set as $\mathcal{X}_i:=\{\mathsf{A}\} \cup \mathcal{M} \cup \{\mathsf{O}\} \cup \mathcal{R}$.
We label as $x_i$ an element of this set.

The elements of the set of commodity labels $\mathcal{C}$ are C, E, and N, which respectively abbreviate corn, ethanol, and natural gas.
The spot price of commodity $c \in \mathcal{C}$ in stage $i$ is $s^c_i \in \mathbb{R}_+$.
The price in stage $i$ of the futures for commodity $c$ with delivery in stage $j \geq i$ is $F^{c}_{i,j} \in \mathbb{R}_+$.
If $i$ equals $j$ then $F^{c}_{i,i}$ and $s^c_i$ coincide.
The stage $i$ forward curve of commodity $c$, ${F}_i^c$, is the vector $(F^{c}_{i,j},j\in\{i,\ldots,I-1\})$.
We define the vector of forward curves in stage $i$ as ${F}_i := ({F}^{c}_i,c\in\mathcal{C})$.
It takes values in $\mathbb{R}^{3(I-i)}_+$.

We denote as $\mathsf{P}$ the decision to produce $Q$ gallons of ethanol taking as inputs $\gamma_{\textrm{C}}$ bushels of corn and $\gamma_{\textrm{N}}$ mmBTU of natural gas each per gallon of output.
We label as $\mathsf{S}$ the production suspension action.
All other choices correspond to modifications of the current operating mode.
Hence, we represent them by employing the notation for the resulting such mode.
For stage $i\in\{0,\ldots,I-2\}$, the feasible action set corresponding to the operating mode $x_i$, $\mathcal{A}_i(x_i)$, is
respectively $\{\mathsf{A},{\mathsf{M}}_1,\mathsf{P},\mathsf{S}\}$ and $\{\mathsf{A},\mathsf{P},\mathsf{S}\}$ when $i\leq I-1 - N^\mathsf{M}$ and $i> I-1 - N^\mathsf{M}$ if $x_i = \mathsf{O}$;
$\{\mathsf{M}_{n+1}\}$ if $x_i =\mathsf{M}_{n}$ and $n \neq N^\mathsf{M}$;
respectively $\{\mathsf{A},{\mathsf{M}}_{N^{\mathsf{M}}},{\mathsf{R}}_{1}\}$ and $\{\mathsf{A},{\mathsf{M}}_{N^{\mathsf{M}}}\}$ when $i\leq I-1 - N^\mathsf{R}$ and $i> I-1 - N^\mathsf{R}$ if $x_i = \mathsf{M}_{N^\mathsf{M}}$;
$\{\mathsf{R}_{n+1}\}$ if $x_i =\mathsf{R}_{n}$ and $n \neq N^\mathsf{R} - 1$;
$\{\mathsf{O}\}$ if $x_i = N^\mathsf{R} - 1$; and, for ease of model formulation,
$\{\mathsf{A}\}$ if $x_i = \mathsf{A}$.
In stage $I-1$, the feasible action set for the operating mode $x_{I-1}$, $\mathcal{A}_{I-1}(x_{I-1})$, includes $\mathsf{A}$ as its only element.

The unitary production spread (gross margin) is $s_i^{\textrm{E}} - \gamma_{\textrm{C}} s_i^{\textrm{C}} - \gamma_{\textrm{N}} s_i^{\textrm{N}}$.
The respective costs per stage of producing $Q$ gallons of ethanol, suspending production, and keeping the plant fully mothballed are ${\mathsf{C}}_{\mathsf{P}}$, ${\mathsf{C}}_{\mathsf{S}}$ ($< {\mathsf{C}}_{\mathsf{P}}$), and ${\mathsf{C}}_{\mathsf{M}}$ ($< {\mathsf{C}}_{\mathsf{S}}$) dollars.
The costs of initiating the mothballing and reactivation processes are respectively $\mathsf{I}_{\mathsf{M}}$ and $\mathsf{I}_{\mathsf{R}}$ dollars.
Abandoning the plant yields a net salvage value of $S$ dollars.
The per stage reward depends on the operating mode $x_i$, the spot price vector $s_i:=(s^c_{i},c\in\mathcal{C})$, and the action $a_i\in\mathcal{A}(x_i)$. 
We define this function as
\begin{align*}
r(x_i,{s}_i,a_i) := \left\{\begin{array}{ll}
(s^{\textrm{E}}_{i} - \gamma_{\textrm{C}} s^{\textrm{C}}_{i} - \gamma_{\textrm{N}} s^{\textrm{N}}_{i})Q - \mathsf{C}_\mathsf{P},&\mbox{if $(x_i,a_i) \in (\mathsf{O},\mathsf{P})$},\\
-\mathsf{C}_{\mathsf{S}},&\mbox{if $(x_i,a_i) = (\mathsf{O},\mathsf{S})$},\\
- \mathsf{I}_{\mathsf{M}},&\mbox{if $(x_i,a_i) = (\mathsf{O},{\mathsf{M}_1})$},\\
- \mathsf{C}_{\mathsf{M}},&\mbox{if $(x_i,a_i) = ({\mathsf{M}_{N^\mathsf{M}}},{\mathsf{M}_{N^\mathsf{M}}})$},\\
-\mathsf{I}_{\mathsf{R}},&\mbox{if $(x_i,a_i) = ({\mathsf{M}_{N^\mathsf{M}}},\mathsf{R}_1)$},\\
S,&\mbox{if $(x_i,a_i)\in \{(\mathsf{O},\mathsf{A}),({\mathsf{M}}_{N^\mathsf{M}},\mathsf{A})\}$},\\
0,&\mbox{if $(x_i,a_i) \in \{(\mathsf{A},\mathsf{A}),({\mathsf{M}}_n,{\mathsf{M}}_{n+1} | n \neq N^\mathsf{M}),$}\\
&\hspace{0.9in}{({\mathsf{R}}_{N^\mathsf{R} - 1},\mathsf{O}), ({\mathsf{R}}_n,{\mathsf{R}}_{n+1} | n \neq N^\mathsf{R} - 1)\}}.\\
\end{array}\right.
\end{align*}

The function $f(x_i,a_i)$ gives the next stage operating mode that results from executing feasible action $a_i$ in the current stage when the operating mode is $x_i$.
In particular, its value is $\mathsf{O}$ if the pair $(x_i,a_i)$ belongs to the set $\{(\mathsf{O},\mathsf{P}),(\mathsf{O},\mathsf{S})\}$ and $a_i$ in all other cases.
The plant operations do not affect the evolution of the vector of forward curves; that is, the plant is small relative to the markets.
The dynamics of this vector are governed by exogenously specified stochastic processes.

The set of feasible policies is $\Pi$.
Such a policy $\pi$ is the collection of decision rules $\{A^{\pi}_i, i\in\mathcal{I}\}$, with $A^{\pi}_i: \mathcal{X}_i \times \mathbb{R}^{3(I-i)}_+ \rightarrow \mathcal{A}_i(x_i)$.
The objective is to choose a feasible policy that maximizes the market value of operating the plant during the finite horizon given the initial state $(x_0,F_0)$:
\begin{equation}
\max_{\pi\in\Pi}\sum_{i\in\mathcal{I}} \delta^i \mathbb{E}\left[r\left(x^\pi_{i},s_{i},A^{\pi}_{i}\left(x^\pi_{i},F_{i}\right)\right)\mid x_0,F_0\right],
\label{eq:MDP}
\end{equation}
where $\delta$ is the per stage risk free discount factor; $\mathbb{E}$ is expectation under a risk neutral measure (see, e.g., \citealp{kn:Shreve2004}) for the stochastic processes that determine the evolution of the vector of forward curves; and $x^\pi_{i}$ is the random operating mode reached in stage $i$ when using policy $\pi$.

\section{Pathwise Optimization}
\label{Section 3}
In this section we present a PO approach to approximately solve MDP \eqref{eq:MDP}. In \S\ref{subsec:OptOfDualPenalties} we develop a dual formulation of this MDP. In \S\ref{subsec:PLP} we introduce PLP to approximate this model. In \S\ref{subsec:estBndsGrdPol} we describe how to use part of the PLP solution to estimate a dual bound on the optimal policy value and obtain a feasible policy, which we use to determine a lower bound on this value.

\subsection{Dual Model}
\label{subsec:OptOfDualPenalties}
We introduce a dual version of MDP \eqref{eq:MDP} in which decisions are made knowing the realized paths of the vectors of forward curves but the advantage of this foresight is entirely eliminated by ideal penalties \citep{brown2010information}.

Let $\bar{F}$ be a sample path that includes the vectors of forward curves from stages $0$ through $I-1$ starting with $F_{0}$ (we suppress this dependence from our notation for ease of exposition). The set $\bar{\mathcal{F}}$ is the collection of all such paths. We denote by $s_i(\bar F)$ and $F_i(\bar{F})$ the stage $i$ vectors of spot prices and forward curves corresponding to sample path $\bar{F}$. The dual policy $\bar \pi$ is the collection of decision rules $\{\bar{A}^{\bar \pi}_i, i \in \mathcal{I}\}$, where $\bar{A}^{\bar \pi}_i: \mathcal{X}_i \times \bar{\mathcal{F}} \rightarrow  \mathcal{A}_i(x_i)$ prescribes a feasible action for stage $i$, operating mode $x_i$, and sample path $\bar{F}$. The set of such policies is $\bar \Pi$.

Ideal penalties depend on the value function associated with model (\ref{eq:MDP}), which solves the stochastic dynamic program \citep{puterman1994markov}
\begin{equation}
V_{i}(x_{i},F_{i}) =\max_{a_{i}\in \mathcal{A}_i({x_{i}})} \left\{r_{i}(x_{i},s_{i},a_{i})+\delta\mathbb{E}\left[V_{i+1}(f(x_{i},a_{i}), F_{i+1})\bigg|F_{i}\right]\right\} \label{SDP}
\end{equation}
for each $(i,x_{i},F_{i})\in \mathcal{I}\setminus\{0\}\times\mathcal{X}_{i}\times\mathcal{F}_{i}$ with boundary conditions $V_I(x_I,F_I) \coloneqq 0$ for each $(x_I, F_I) \in \mathcal{X}_I \times \mathcal{F}_I$. Consider stage $i\neq I-1$ and suppose we take feasible action $a_{i}$ for operating mode $x_{i}$ and sample path $\bar{F}$. The ideal penalty is the additional value of knowing the stage $i+1$ forward curve $F_{i+1}(\bar F)$ at stage $i$ relative to only having knowledge of the forward curve $F_i(\bar{F})$ at this stage, which corresponds to the discounted difference
\begin{equation}
\delta \left(V_{i+1}\left(f(x_i,a_{i}),F_{i+1}(\bar{F})\right) - \mathbb{E}[V_{i+1}\left(f(x_i,a_{i}),F_{i+1}\right)|F_i(\bar{F})]\right).\label{DualPenalty}
\end{equation}

We use the penalty (\ref{DualPenalty}) to reduce the cash flow that ensues the stage $i\neq I-1$ from applying the decision rule $\bar{A}^{\bar \pi}_{i}$ to the pair $(x_i,\bar F)$. The resulting dual MDP is\looseness=-1
\begin{align}
&\mathbb{E}\left[\max_{\bar \pi\in \bar{\Pi}} \bigg\{\sum_{i\in\mathcal{I}\setminus\{I-1\}}\delta^{i}\bigg[r(x^{\bar \pi}_{i},s_{i}(\bar{F}),\bar{A}^{\bar \pi}_{i}) - \delta \big(V_{i+1}(f(x^{\bar \pi}_i,\bar{A}^{\bar\pi}_i),F_{i+1}(\bar{F}))\right.\nonumber\\
&\hspace{2.9in}\left.  - \mathbb{E}[V_{i+1}(f(x^{\bar \pi}_i,\bar{A}^{\bar \pi}_i),F_{i+1})|F_i(\bar{F})]\big) \bigg]\right. \nonumber\\
& \hspace{2.9in}\left. +\delta^{I-1}r(x^{\bar\pi}_{I-1},s_{I-1}(\bar{F}),\bar{A}^{\bar\pi}_{I-1})\bigg\}\bigg|x_{0},F_{0}\right],
\label{dualMDP}
\end{align}
where we use the shorthand notation $\bar{A}^{\bar\pi}_{i}$ instead of $\bar{A}^{\bar\pi}_{i}(x^{\bar\pi}_i,\bar{F}_i)$. This model differs from MDP \eqref{eq:MDP} in two key ways: (i) The maximization is inside the expectation because dual policies depend on sample paths and (ii) its objective function is the sum of the discounted ideally penalized rewards and the last stage reward. Let $V_0(x_0,F_0)$ be the value function for stage 0 and the given state $(x_{0},F_{0})$, which is obtained in a manner analogous to \eqref{SDP} for this stage and state. At optimality the objective function (\ref{dualMDP}) equals $V_0(x_0,F_0)$ for each sample path \citep{brown2010information}. It follows that the MDP \eqref{eq:MDP} and its version \eqref{dualMDP} are equivalent. \looseness=-1

\subsection{Pathwise Linear Program}\label{subsec:PLP}
The dual model \eqref{dualMDP} is intractable because (i) the outer expectation is impossible to evaluate exactly in general, and in particular for our application discussed in \S\ref{Section 5}, and (ii) the ideal penalties are unknown. We formulate PLP to address both these issues.  

First, we approximate the outer expectation with a sample average based on $L$ Monte Carlo simulation sample paths of the vectors of forward curves from stage $0$ through stage $I-1$ starting from $F_{0}$. We define the index set $\mathcal{L}\coloneqq\{1,...,L\}$. We denote by $F^l_i$ and $s^l_i$, respectively, the stage $i$ vectors of forward curves and spot prices for sample path $l\in\mathcal{L}$. Second, we replace the ideal penalties in \eqref{dualMDP} by ``good'' penalties based on VFAs  \citep{brown2010information}. We do so in a parsimonious manner that avoids the specification of a VFA, and thus penalties, at states with fixed actions, that is, both when the plant is abandoned and when it is being mothballed or reactivated. In other words, we only need to approximate the value function for operating modes $\mathsf{M}_{N^\mathsf{M}}$ and $\mathsf{O}$, which we include in set $\mathcal{X}'$. We specify a VFA for stage $i \in \mathcal{I}\setminus\{0, I-1\}$ and state $(x_{i},F_i) \in \mathcal{X}' \times \mathcal{F}_{i}$ as
$\sum_{b\in\mathcal{B}_{i}}\beta_{i, x_{i},b}\phi_{i,b}(F_{i})$,
where $\mathcal{B}_{i}:=\{1,...,B_i\}$ is the index set of $B_i$ basis functions for stage $i$, $\beta_{i, x_{i}, b}\in\mathbb{R}$ is the weight associated with the $b$-th basis function for stage $i$ when the operating mode is $x_{i}$, and $\phi_{i,b}(F_i)$ is this function whose argument is the vector of forward curves $F_i$ . We introduce the modified transition function $f'(x_i,a_i)$, which equals $\mathsf{O}$ if   $(x_{i},a_{i})\in\{(\mathsf{O},\mathsf{P}), (\mathsf{O},\mathsf{S}), (\mathsf{M}_{N^{\mathsf{M}}},\mathsf{R_1})\}$, $\mathsf{M}_{N^{\mathsf{M}}}$ if $(x_{i},a_{i}) \in \{(\mathsf{O},\mathsf{M}_1),(\mathsf{M}_{N^{\mathsf{M}}}, \mathsf{M}_{N^{\mathsf{M}}})\}$, and $\mathsf{A}$ if $(x_{i},a_{i})\in \{(\mathsf{O},\mathsf{A}),(\mathsf{M}_{N^{\mathsf{M}}},\mathsf{A}), (\mathsf{A},\mathsf{A})\}$. Further, we define the stage transition function $g(i,x_i,a_i)$, with $x_{i}\in\mathcal{X}'$, to be equal to $i + N^\mathsf{M}$ if $(x_i, a_i) = (\mathsf{O}, \mathsf{M}_1)$, $i + N^\mathsf{R}$ if $(x_i, a_i) = (\mathsf{M}_{N^{\mathsf{M}}},\mathsf{R}_1$), and $i+1$ otherwise. We replace the value function with our VFA in the ideal dual penalty (\ref{DualPenalty}) to obtain the following good dual penalty for stage $i\in\mathcal{I}\setminus\{0, I-1\}$, operating mode $x_i \in \mathcal{X}'$, sample path $l\in\mathcal{L}$, and action $a_{i}\in \mathcal{A}_{i}(x_i)$:
\begin{equation}\sum_{b\in\mathcal{B}_i}\beta_{g(i,x_i,a_i), f'(x_{i},a_i),b}\Delta_i^{\mathbb{E},l} \phi_{g(i,x_i,a_i),b},\label{feasDualPenalty}
\end{equation}
where $\Delta_i^{\mathbb{E},l} \phi_{g(i,x_i,a_i),b} := \delta^{g(i,x_i,a_i)-i}\{\phi_{g(i,x_i,a_i),b}(F^l_{g(i,x_i,a_i)})-\mathbb{E}[\phi_{g(i,x_i,a_i),b}(F_{g(i,x_i,a_i)})|F^l_{i}]\}.$ For a fixed VFA weight vector $\beta$, let $U_{0}^{l,\beta}(x_{0})$, with $x_{0}\in\mathcal{X}'$, be the optimal objective function value of the following maximization for sample path $l\in\mathcal{L}$:
\begin{align}
\max_{x^l,a^l}&\sum_{i\in\mathcal{I}\setminus\{I-1\}}\delta^{i}\bigg[r(x^l_{i},s^l_{i},a^l_{i}) - \mathbbm{1}(f'(x_{i}^{l},a_{i}^{l})\in\mathcal{X}')\sum_{b\in\mathcal{B}_{i}}\beta_{g(i,x^l_i,a^l_i), f'(x^l_{i},a^l_i),b}\Delta_i^{\mathbb{E},l} \phi_{g(i,x^l_i,a^l_i),b}\bigg]\nonumber\\
&\ \ \ +\delta^{I-1}r(x^l_{I-1},s^l_{I-1},a^l_{I-1}),\label{POOptModel}\\
\mbox{s.t.}&\ \ \  x_{0}^{l}=x_{0}, \label{POOptModel: InitialState}\\
&\ \ \  x_{i+1}^{l}=f'(x_{i}^{l},a_{i}^{l}), \forall i\in\mathcal{I}\setminus\{I-1\};\label{POOptModel: StateTrans}
\end{align}
where $\mathbbm{1}(\cdot)$ is the indicator function that evaluates to $1$ when its argument is true and to $0$ otherwise and
$x^l$ and $a^l$ are the vectors of operating modes and actions, respectively, for sample path $l$. The dual model resulting from our approximations is \looseness=-1
\begin{align}
\label{POOptiModel:OBJ}
\min_{\beta}\frac{1}{L}\sum_{l\in\mathcal{L}}U^{l,\beta}_{0}(x_0).
\end{align}
Each term $U_{0}^{l,\beta}(x_{0})$ can be obtained by first solving the dynamic program
\begin{align}
{U}_{i}^{l,\beta}(x_{i}) &= \max_{a_i \in \mathcal{A}_i(x_i)}\bigg\{ r(x_i,s_i^{l},a_i) -  \mathbbm{1}(f'(x_i,a_i)\in\mathcal{X}')\sum_{b\in\mathcal{B}_i}\beta_{g(i,x_i,a_i), f'(x_{i},a_i),b}\Delta_i^{\mathbb{E},l} \phi_{g(i,x_i,a_i),b}\nonumber\\
&\hspace{0.8in} + \delta {U}_{g(i,x_i,a_i)}^{l,\beta}(f'(x_{i},a_i))\bigg\},\label{eq:samplePathDualDP2}
\end{align}
for each pair $(i,x_i) \in \mathcal{I} \setminus \{0,I-1\} \times \mathcal{X}^{''}$ where $\mathcal{X}^{''}:=\mathcal{X}'\cup\{\mathsf{A}\}$, with boundary conditions
\begin{equation}
{U}_{I - 1}^{l,\beta}(x_{I-1}) = \max_{a_{I-1} \in \mathcal{A}_{I-1}(x_{I-1})}r(x_{I-1},s^{l}_{I-1},a_{I-1}),\label{eq:samplePathDualDP3}
\end{equation}
 for each $x_{I-1} \in \mathcal{X}^{''}$, and then performing the analogous optimization on the right hand side of \eqref{eq:samplePathDualDP2} for stage $0$ and the given starting operating mode $x_{0}$. Each term $U_{0}^{l,\beta}(x_{0})$ and the value function of each such dynamic program solve the following linear program \citep{manne1960linear}: 
 \begin{align}
\min_{{U}^l} \ \ {U}_{0}^{l}(x_{0})&\label{DPLPObj}\\
\mbox{s.t.} \ \
 {U}_{0}^{l}(x_{0})&\geq r(x_{0},s_{0}^{l},a_{0}) - \mathbbm{1}(f'(x_0,a_0)\in\mathcal{X}')\sum_{b\in\mathcal{B}_i}\beta_{g(0,x_0,a_0), f'(x_{0},a_0),b}\Delta_0^{\mathbb{E},l} \phi_{g(0,x_0,a_0),b} \notag\\
&\hspace{0.2in}+\delta {U}_{g(0,x_0,a_0)}^{l}(f'(x_{0},a_{0})), \forall a_{0} \in \mathcal{A}_{0}(x_{0}),\label{DPLPConst1}\\
 {U}_{i}^{l}(x_{i})&\geq r(x_{i},s_{i}^{l},a_{i})- \mathbbm{1}(f'(x_i, a_i)\in\mathcal{X}')\sum_{b\in\mathcal{B}_i}\beta_{g(i,x_i,a_i), f'(x_{i},a_i),b}\Delta_i^{\mathbb{E},l} \phi_{g(i,x_i,a_i),b}\nonumber\\
&\hspace{0.2in}+\delta {U}_{g(i,x_i,a_i)}^{l}(f'(x_{i},a_{i})), \forall (i,x_i,a_i) \in \mathcal{I} \setminus \{0,I-1\} \times \mathcal{X}'' \times \mathcal{A}_i(x_i),\label{DPLPConst2}\\
{U}_{I-1}^{l}(x_{I-1})&\geq r(x_{I-1},s_{I-1}^{l},a_{I-1}), \ \ \ \forall (x_{I-1},a_{I-1}) \in \times \mathcal{X}'' \times \mathcal{A}_{I-1}(x_{I-1}),\label{DPLPConst3}
\end{align}
where the vector $U^{l}$ includes the variables $U_{i}^{l}(x_{i})$ for the pair $(0,x_0)$ and each triple $(i,x_i,l)\in\mathcal{I}\setminus\{0\}\times\mathcal{X}''\times\mathcal{L}$. This model minimizes the value of the variable $U^l_0(x_0)$, which is the objective function \eqref{DPLPObj}, by imposing on its decision variables relaxed versions of both the analogue of equation (\ref{eq:samplePathDualDP2}) for the pair $(0,x_0)$ and the conditions \eqref{eq:samplePathDualDP2}-\eqref{eq:samplePathDualDP3}, that is, the sets of constraints \eqref{DPLPConst1}-\eqref{DPLPConst3}.  	Replacing the quantity $U^{l,\beta}_{0}(x_0)$ in \eqref{POOptiModel:OBJ} with the minimization \eqref{DPLPObj}-\eqref{DPLPConst3} for each sample path $l \in \mathcal{L}$ allows us to equivalently express \eqref{POOptiModel:OBJ} as the following linear program, which is PLP:	
\begin{align}
\min_{\beta,{U}} \ \ \frac{1}{L}\sum_{l\in\mathcal{L}}{U}_{0}^{l}(x_{0})&\label{PLPObj}\\
\mbox{s.t.} \ \
 {U}_{0}^{l}(x_{0})&\geq r(x_{0},s_{0}^{l},a_{0}) - \mathbbm{1}(f'(x_0,a_0)\in\mathcal{X}')\sum_{b\in\mathcal{B}_i}\beta_{g(0,x_0,a_0), f'(x_{0},a_0),b}\Delta_0^{\mathbb{E},l} \phi_{g(0,x_0,a_0),b} \notag\\
&\hspace{0.2in}+\delta {U}_{g(0,x_0,a_0)}^{l}(f'(x_{0},a_{0})), \forall (l,a_{0}) \in \mathcal{L}\times \mathcal{A}_{0}(x_{0}),\label{PLPConst1}\\
 {U}_{i}^{l}(x_{i})&\geq r(x_{i},s_{i}^{l},a_{i})- \mathbbm{1}(f'(x_i,a_i)\in\mathcal{X}')\sum_{b\in\mathcal{B}'}\beta_{g(i,x_i,a_i), f'(x_{i},a_i),b}\Delta_i^{\mathbb{E},l} \phi_{g(i,x_i,a_i),b}\nonumber\\
&\hspace{0.2in}+\delta {U}_{g(i,x_i,a_i)}^{l}(f'(x_{i},a_{i})), \forall (l,i,x_i,a_i) \in \mathcal{L}\times \mathcal{I} \setminus \{0,I-1\} \times \mathcal{X}^{''}\times \mathcal{A}_i(x_i),\label{PLPConst2}\\
 {U}_{I-1}^{l}(x_{I-1})&\geq r(x_{I-1},s_{I-1}^{l},a_{I-1}), \ \ \ \forall (l,x_{I-1},a_{I-1}) \in \mathcal{L}\times \mathcal{X}^{''} \times \mathcal{A}_{I-1}(x_{I-1}).\label{PLPConst3}
\end{align}
Proposition \ref{PLPProp} establishes that PLP is well defined.  
\begin{proposition}
\label{PLPProp}
PLP has a finite optimal objective function value and at least one bounded optimal solution.
\end{proposition}

\subsection{Dual Bound, Greedy Policy, and Lower Bound}
\label{subsec:estBndsGrdPol}
	The optimal PLP objective function value is not a valid dual bound in general because it suffers from a sample average optimization bias. We thus estimate an unbiased dual bound using a new set of Monte Carlo simulation sample paths of the vectors of forward curves and the VFA coefficient vector $\beta^{\mathrm{PLP}}$ obtained by solving PLP. This estimation can be performed by obtaining the analogue of $U_{0}^{l,\beta}(x_{0})$ based on the dynamic program \eqref{eq:samplePathDualDP2}-\eqref{eq:samplePathDualDP3} but using $\beta^{\mathrm{PLP}}$ and these sample paths in lieu of $\beta$ and the ones indexed by the elements of set $\mathcal{L}$, respectively, and averaging each such quantity. The dual penalty terms that appear in these dynamic programs include expectations that need to be evaluated. Approximating them by sample average approximations is a possibility \citep{desai2012pathwise} but introduces an error in the dual bound estimate. We thus choose basis functions and stochastic models for the evolution of the vector of forward curves that satisfy Assumption \ref{ass:closedform}, which is common in the literature (see, e.g., \citealt{nadarajah2017comparison} and references therein)
\begin{assumption}
\label{ass:closedform}
The expectation $\mathbb{E}[\phi_{j,b}(F_{j})|F_i]$ is available in an efficiently computable closed form for each $i$ and $j \in \mathcal{I}\setminus\{I-1\}$ with $j > i$ and $F_i \in \mathcal{F}_i$.
\end{assumption}
To obtain an operating policy and estimate a lower bound, we employ the well-known ``greedy'' optimization framework in approximate dynamic programming (see, e.g., \citealt{kn:Powell2011}, \S 6.4). Given a VFA weight vector $\beta$, the stage $i$ greedy decision rule for states that include operating modes that admit more than one feasible action (that is, $x_i \in \mathcal{X}'$) is\looseness=-1
\begin{equation*}
\argmax_{a_{i} \in \mathcal{A}_i(x_i)}\left\{r(x_{i},s_{i},a_{i})+\delta \sum_{b\in\mathcal{B}_i}\beta_{g(i,x_i,a_i), f'(x_{i},a_i),b}\mathbb{E}[\phi_{g(i,x_i,a_i),b}(F_{g(i,x_i,a_i)})|F_i]\right\},
\end{equation*} 
with ties broken in some prespecified way. The collection of these decision rules is the greedy policy. To estimate its associated lower bound, we employ the same set of sample paths of the vectors of forward curves used to obtain an unbiased dual bound estimate and apply the greedy policy to the states visited along each such path starting from the initial stage and state. The average of the sum of the resulting discounted rewards is an unbiased lower bound estimate.\looseness=-1

Although a VFA coefficient vector obtained by solving PLP can be used to derive a greedy policy, its corresponding lower bound estimate may be weak \citep{desai2012pathwise}. Indeed, consider the common assumption that the first basis function used to construct the VFAs is a constant, that is, $\phi_{i,1}(\cdot) \coloneqq 1$ for each $i \in \mathcal{I}\setminus\{0,I-1\}$, which implies that each term $\Delta_i^{\mathbb{E},l} \phi_{g(i,x_i,a_i),1}$ equals zero. Thus, the decision variables $\beta_{i,x_i, 1}$'s have zero coefficients in PLP and the resulting VFAs do not have intercepts, which is undesirable from the perspective of obtaining a good greedy policy. To address this issue, following \citet{desai2012pathwise}, we determine VFAs with intercepts based on the vector $U^{\mathrm{PLP}}$ obtained by solving PLP. Let $U_{i}^{l,\beta^{\mathrm{PLP}}}(x_{i},F_{i}^{l})$ be an element of $U^{\mathrm{PLP}}$. For each pair $(i,x_i) \in \mathcal{I}\setminus\{0,I-1\} \times \mathcal{X}'$, we define the vector $\beta_{i,x_i} := (\beta_{i,x_i,b},b\in\mathcal{B}_i)$ and solve the regression model
\begin{align*}
\min_{\beta_{i,x_{i}}}\frac{1}{L}\sum_{l\in\mathcal{L}}\left(U_{i}^{l,\beta^{\mathrm{PLP}}}(x_{i})-\sum_{b\in\mathcal{B}_i}\beta_{i,x_{i},b}\phi_{i,b}(F_i^l)\right)^{2}.
\end{align*}
We employ these resulting optimal solutions to specify VFAs and consequently obtain a greedy policy, from which we estimate a lower bound.
\section{Solving the Pathwise Linear Program}
\label{Section 4}
PLP is both ill conditioned and large scale in our application. We present the pre-conditioning and optimization algorithms that we develop to address these issues in \S\ref{sec:BPCA} and in \S\ref{subsec:BCD}, respectively.\looseness=-1

\subsection{Pre-conditioning Algorithm}\label{sec:BPCA}
The PLP columns that correspond to the coefficients of the basis functions are nearly parallel in our numerical study, which is based on commonly used such functions (see, e.g., \citealp{longstaff2001valuing,boogert2011gas, nadarajah2017comparison}, and references therein). That is, the resulting PLP is ill-conditioned and thus difficult to solve. For example, Gurobi, a state-of-the-art commercial optimization solver, faces severe numerical issues and is unable to find an optimal solution when applied to even medium sized PLP instances. In contrast, Gurobi readily solves these instances to optimality after we execute on them the (computationally efficient) procedure discussed below.\looseness=-1

To simplify exposition, we express the PLP constraints in matrix form as
\begin{equation}
Q{U} + G \beta \geq r, \label{PLPform:Const}
\end{equation} 
where ${U}$ and $r$ are the column vectors $({U}^l_i(x_i),(l,i,x_i) \in\mathcal{L}\times\mathcal{I} \times \mathcal{X}'')$ and $(r(x_i,s_i^l,a_i), (l,i,x_i,a_i) \in \mathcal{L}\times \mathcal{I} \times \mathcal{X}'' \times \mathcal{A}(x_i))$, respectively, and $Q$ and $G$ are the respective constraint coefficient matrices associated with the decision variable vectors ${U}$ and $\beta$. The row of the $Q$ matrix indexed by $(l, i, x_i, a_i)$ has nonzero coefficients (respectively equal to one and $-\delta$) only in the columns for the variables ${U}^l_i(x_i)$ and ${U}_{g(i,x_i,a_i)}^{l}(f'(x_{i},a_{i}))$. The $G$ matrix has a block-diagonal structure. For each pair $(i,x_i)\in\mathcal{I}\setminus\{0,I-1\}\times\mathcal{X}'$, we define the set of tuples $\mathcal{T}_i(x_i) := \{(l,j,x_{j},a_{j}) | (l,j,x_{j},a_{j}) \in\mathcal{L}\times \mathcal{I}\setminus \{i, I-1\} \times \mathcal{X}' \times \mathcal{A}(x_j), g(j,x_j,a_j) = i, f'(x_{j},a_j) = x_i\}$ to describe it. The $(i,x_i)$-th block of $G$, $G_{i,x_i}$, includes the columns associated with the triples $(i,x_i,b)$'s for each $b\in\mathcal{B}_i$ and the rows corresponding to the tuples $(l,j,x_{j},a_{j})$'s in set $\mathcal{T}_{i}(x_i)$. Figure \ref{fig:MatrixStructure} illustrates this structure. \looseness=-1
\begin{figure}[h!]
	\[ \left[ \begin{array}{cccc}
	G_{1, \mathsf{O}} & 0 & \cdots & 0 \\
	0 & G_{1, \mathsf{M}_{N^{\mathsf{M}}}} & \cdots & 0 \\
	\vdots  & \vdots & \ddots & \vdots \\
	0 & 0 & \cdots & G_{I-2, \mathsf{M}_{N^{\mathsf{M}}}} \end{array} \right]\] 	
\caption{Structure of the G matrix}
\label{fig:MatrixStructure}
\end{figure}

PLP is ill-conditioned because the columns of the $G$ matrix are almost parallel. We use PCA to make the columns of each block $G_{i,x_i}$ of this matrix perpendicular. We refer to this procedure as block PCA (BPCA). It both preserves the block-diagonal structure of the $G$ matrix and has a smaller computational burden than applying PCA directly to this matrix. Specifically, we denote by $W_{i,x_i}$ the square matrix with columns equal to the eigenvectors of $G_{i,x_i}^{\mathsf{T}}G_{i,x_i}$, where the superscripted $\mathsf{T}$ denotes transposition, and use it to obtain $G^{\scriptscriptstyle \perp}_{i,x_i}$ as the orthogonal linear transformation $G_{i,x_i}W_{i,x_i}$ of $G_{i,x_i}$. 
We denote by $G^{\scriptscriptstyle \perp}$ the analogue of $G$ with each block $G_{i,x_i}$ replaced by $G^{\scriptscriptstyle \perp}_{i,x_i}$. We then replace the PLP constraints (\ref{PLPform:Const}) with 
\begin{equation}
Q{U} + G^{\scriptscriptstyle \perp} \beta \geq r \label{form:Const}
\end{equation}
We refer to the resulting linear program as the preconditioned PLP (P2LP). Proposition \ref{prop:equivPCAPLP} establishes the equivalence of PLP and P2LP at optimality.
\begin{proposition}
\label{prop:equivPCAPLP}
PLP and P2LP have identical sets of optimal solutions.
\end{proposition}
\subsection{Optimization Algorithm}\label{subsec:BCD}
P2LP is a large scale model in our application. Specifically, it has three million variables and ten million constraints. Attempting to solve such P2LP instances using a commercial solver requires too much memory in our numerical study. We thus devise a customized solution method to deal with this issue.

Conceptually, our approach is a cyclic BCD (CBCD) algorithm that aims at solving (\ref{POOptiModel:OBJ}).
This model corresponds to the unconstrainted minimization of a piecewise linear convex function that is not explicitly available. That is, evaluating this function requires solving (\ref{DPLPObj}) - (\ref{DPLPConst3}) for each sample path of the vector of forward curves $l\in\mathcal{L}$. Thus, our technique cyclically optimizes a sequence of P2LPs in which the values of some of the decision variables that belong to the $\beta$ vector are fixed and the ones of all the others, including those that are part of the $U$ vector, are optimally chosen. In particular, exact evaluation of the objective function (\ref{POOptiModel:OBJ}) requires an optimal selection of the values of all elements of the $U$ vector.



\begin{algorithm}[h!]
\caption{CBCD Algorithm}
\belowdisplayskip=0pt
\label{Alg:BCDPO}
\SetKwRepeat{Do}{do}{while}
\SetKwInOut{Input}{input}\SetKwInOut{Output}{output}
\SetKwInOut{Initialization}{initialization}
\Input{Initial vector $\beta^0$, $\mathrm{OBJ}(\beta^0)$, block partition $\mathcal{P}$, and stopping tolerance $\epsilon > 0$.}
\Initialization{Set $h = 0$.}
\Do{$|\mathrm{OBJ}(\beta^h) - \mathrm{OBJ}(\beta^{h-1})| > \epsilon$}{
$h = h + 1$.\\
\For{$p = 1$ to $P$}{
(i) Let $\beta^{h,p}$ be an optimal $\beta$ vector for the linear program\vspace{-0.2in}
\begin{align}
\min_{\beta, U}&\frac{1}{L}\sum_{l\in\mathcal{L}}{U}_{0}^{l}(x_{0})\label{form:ObjBCDLP}\\
\mbox{s.t. } &Q {U} + G^{\scriptscriptstyle \perp} \beta \geq r,\label{form:ConstBCDLP1}\\
&\beta(\mathcal{P}_{p'}) = \beta^{h}(\mathcal{P}_{p'}), \forall p' \in \{1,\ldots,p-1\},\\
&\beta(\mathcal{P}_{p'}) = \beta^{h-1}(\mathcal{P}_{p'}), \forall p' \in \{p+1,\ldots,P\}.\label{form:ConstBCDLP2}
\end{align} 
(ii) $\beta^{h}(\mathcal{P}_{p}) := \beta^{h,p}$.
}
\vspace{0.1in}
Evaluate $\mathrm{OBJ}(\beta^h)$.
}
\Output{Return $\beta^{h}$}
\end{algorithm} 
Our CBCD procedure uses blocks of $\beta$ vector variables that correspond to a partition $\mathcal{P}$ of the set $\mathcal{I}\setminus\{0,I-1\} \times \mathcal{X}'$ into $P\leq\sum_{i\in\mathcal{I}\setminus\{0,I-1\}}|\mathcal{X}'|$ sets $\mathcal{P}_1$ through $\mathcal{P}_{P}$. We define $\beta(\mathcal{P}_{p}) := (\beta_{i,x_i,b}, (i,x_i) \in \mathcal{P}_p, b\in\mathcal{B}_{i})$ as the $p$-th such block. Algorithm \ref{Alg:BCDPO} summarizes our CBCD method. Its inputs are the initial decision variable vector $\beta^0$, $\mathrm{OBJ}(\beta^0)$, where $\mathrm{OBJ}(\cdot)$ is the value of the objective function of (\ref{POOptiModel:OBJ}), the partition $\mathcal{P}$ of the index set of this vector, and a stopping tolerance $\epsilon$. The initialization step sets the iteration counter $h$ to zero. In each subsequent iteration $h$, Algorithm \ref{Alg:BCDPO} (i) selects each set $\mathcal{P}_p$ of the partition $\mathcal{P}$ and solves the variant \eqref{form:ObjBCDLP}-\eqref{form:ConstBCDLP2} of P2LP in which the values of the variables that belong to the $U$ and $\beta(\mathcal{P}_p)$ vectors are optimized, whereas the ones of the others are fixed given the progress of the procedure up to that point and (ii) makes $\beta^{h}(\mathcal{P}_p)$ equal to its corresponding part of the optimal solution to this linear program, $\beta^{h,p}$. Once the solution $\beta^h$ is available, its value $\mathrm{OBJ} (\beta^h)$ is determined. Termination occurs if and only if the quantities $\mathrm{OBJ} (\beta^h)$ and $\mathrm{OBJ} (\beta^{h-1})$ differ by less than $\epsilon$, in which case the vector $\beta^h$ is returned.\looseness=-1

The current literature on BCD algorithms with a cyclic order (see, e.g., \citealp{bertsekas2015convex}, \S6.5) assumes a bounded domain for the decision variables and a unique optimal solution for the model solved in each iteration to prove that they reach an optimal limit point. We establish an analogous result for an idealized version of our CBCD method that satisfies similar assumptions. Suppose there exist (i) finite constants $\beta_{i,x_i,b}^{\mathrm{L}}$ and $\beta_{i,x_i,b}^{\mathrm{H}}$, where the superscripted L and H respectively abbreviate low and high, and (ii) a P2LP optimal solution with $\beta$ vector component $\beta^{*}$ that satisfies the inequalities $\beta_{i,x_i,b}^{\mathrm{L}} \leq \beta^*_{i,x_i,b} \leq \beta_{i,x_i,b}^{\mathrm{H}}$ for each triple $(i,x_i,b) \in \mathcal{I}\setminus\{0, I-1\} \times\mathcal{X}' \times \mathcal{B}_i$ (see Propositions \ref{PLPProp} and \ref{prop:equivPCAPLP}). The idealized CBCD procedure uses zero stopping tolerance and solves the idealized version of P2LP that includes the inequalities 
\begin{equation}\label{betaBounds}
\beta_{i,x_i,b}^{\mathrm{L}} \leq \beta_{i,x_i,b} \leq \beta_{i,x_i,b}^{\mathrm{H}}, \quad \forall (i,x_i,b) \in \mathcal{I}\setminus\{0,I-1\} \times\mathcal{X}'\times\mathcal{B}_i.
\end{equation}
Proposition \ref{PBCDConverges} characterizes the behavior of this method. 
\begin{proposition}
\label{PBCDConverges}
The idealized CBCD algorithm converges to a limit vector. If this vector strictly satisfies the constraints \eqref{betaBounds} and each of its corresponding optimal solutions to (\ref{DPLPObj})-(\ref{DPLPConst3}) are non-degenerate then this vector and the collection of all such solutions optimally solve P2LP.
\end{proposition}
In our numerical study we use the CBCD procedure rather than its idealized version. We find that it always terminates with solutions of seemingly high quality, that is, they lead to greedy policies with small estimated optimality gaps. 


\section{Numerical Study}
\label{Section 5}
In this section we numerically evaluate the performance of the PO approach on the ethanol production application presented in \S\ref{Sec:Model}. We describe the instances used in \S\ref{Section 5.1} and discuss our results in \S\ref{Section 5.2}.
\subsection{Instances}
\label{Section 5.1}
We employ the twelve instances of \cite{nadarajah2018least}. Table \ref{tab:ops_parameters} reports the values of the parameters that are common to all instances. In particular, each instance has twenty four monthly stages. The plant is initially operational ($x_0$ equals $\mathsf{O}$). The starting date of each instance is the first day of each month in 2011. Each initial vector of forward curves ($F_0$) is based on the prices of corn, ethanol, and natural gas New York Mercantile Exchange (NYMEX) futures observed on each such date. The monthly risk-free discount factor ($\delta$) for each beginning date is derived from the one year United States Treasury rate observed on this date. We refer to each instance using the first three letters of the month of their respective initial date.  

\begin{table}[]
    \centering
     \caption{Values of the common parameters.}
    \label{tab:ops_parameters}
\begin{tabular}{cccc}
\hline
     Parameter& Value&Parameter& Value \\
     \hline
     $I $&24 \text{ months}&$I_{\mathsf{M}}$&0.5\\
     $\gamma_{\mathsf{C}}$&0.36 \text{MMBtu/bushel}&$I_{\mathsf{R}}$&2.5\\
     $\gamma_{\mathsf{N}}$&0.035 \text{ MMBtu/gallon}&$\mathsf{C}_{\mathsf{P}}$&2.25\\
     $N^{\mathsf{M}}$&1\text{ month}&$\mathsf{C}_{\mathsf{S}}$&0.5208\\
     $N^{\mathsf{R}}$&3\text{ months}&$\mathsf{C}_{\mathsf{M}}$&0.02917\\
     $Q$ &8.33 \text{million gallons}&${\mathsf{S}}$&0\\
     \hline
\end{tabular}
\end{table}


A typical continuous time stochastic model (\citealp{cortazar1994valuation,clewlow2000energy,blanco2002multi}, \citealp[Chapter 4]{secomandi2014real}) governs the evolution of the vector corn, ethanol, and natural gas forward curves. Denote by $T_i$ the date associated with stage $i \in \mathcal{I}$.  For this model, let $F^{c}(t,T_i)$ be the price of the futures for commodity $c \in \mathcal{C}$ at time $t \in [T_0,T_i]$ with maturity on date $T_i\geq t$. A set of $K$ common factors drives the dynamics of these prices. The $k$-th factor is the standard Brownian motion increment $dW_{k}(t)$. These increments are uncorrelated, that is, $dW_{k}(t)dW_{k^{'}}(t)=0$ for $k,k^{'}\in\mathcal{K}:= \{1,2,...,K\}$ with $k\neq k^{'}$. The time $t$ loading coefficient on the $k$-th factor for the price of the commodity $c$ futures with delivery at time $T$ is $\sigma_{c,k}(t,T)$. The stochastic differential equation that defines the dynamics of the forward curves of commodity $c$ is 
\begin{equation*}
\label{pricemodel}
\frac{dF^{c}(t,T)}{F^{c}(t,T)}=\sum_{k\in\mathcal{K}}\sigma_{c,k}(t, T)dW_{k}(t).
\end{equation*}
This model is dirftless because it is specified under a risk-neutral probability measure. For $i,j \in \mathcal{I}$ with $j > i$, the futures price $F^c_{i,j}$ and spot price $s^c_{i}$ correspond to $F^c(T_i,T_j)$ and $F^c(T_i,T_i)$, respectively. We use the loading coefficient estimates of \cite{nadarajah2018least}, which are based on NYMEX data.   
\subsection{Results}
\label{Section 5.2}

We implement the CBCD algorithm using the same basis functions employed by \cite{nadarajah2018least}. Define $\mathcal{I}_i := \{i,i+1,\ldots,I-1\}$. For each stage $i\in\mathcal{I}$ these functions are (i) one; (ii) $\{F_{i,j}^{c}, j\in \mathcal{I}_{i}, c\in \mathcal{C}\}$; (iii) $\{(F_{i,j}^{c})^2, j\in \mathcal{I}_{i}, c\in \mathcal{C}\}$; (iv) $\{F_{i,j}^{c}F_{i,j}^{c'}, j\in \mathcal{I}_{i}, c, c'\in \mathcal{C}, c\neq c'\}$; and (v) $\{F_{i,j}^{c}F_{i,j+1}^{c}, j\in \mathcal{I}_{i}\setminus\{I-1\}, c\in \mathcal{C}\}$. The conditional expectations of the basis functions of the vector of forward curves are available in \cite{nadarajah2018least}.
They satisfy Assumption \ref{ass:closedform}. We partition the $\beta$ vector of decision variables into four blocks that contain the VFA coefficients for stages 0-5, 6-11, 12-17, and 18-23, respectively. We employ a number of vector of forward curves sample paths ($L$) equal to 70,000: It is the largest value of this parameter that yields P2LP formulations to which we can apply the CBCD algorithm on our high performance computer, discussed below, without facing memory issues. In particular, the ensuing P2LP models have about three million variables and ten million constraints. We use a value for the termination parameter ($\epsilon$) equal to $10^{-3}$. \looseness=-1
\begingroup
\setlength{\tabcolsep}{10pt} 
\renewcommand{\arraystretch}{0.6} 
\begin{table}[h!]
	\centering
	\caption{LSM- and PO-based dual bound estimates, with standard errors reported in parenthesis, and their percentage ratios (the latter ones divided by the former ones, respectively).}
	\label{Table 1}
	\begin{tabular}{cccccc}
		\hline
		&        &                         &                        & & Percentage\\
		Instance & \multicolumn{2}{c}{LSM} & \multicolumn{2}{c}{PO} & Ratio \\ \hline 
		Jan      & 20.96      & (0.03)     & 19.90     & (0.08)     & 95     \\
		Feb      & 20.22      & (0.03)     & 19.32     & (0.08)     & 96     \\
		Mar      & 25.17      & (0.03)     & 23.96     & (0.08)     & 95     \\
		Apr      & 26.64      & (0.03)     & 25.45     & (0.08)     & 96     \\
		May      & 23.01      & (0.03)     & 21.82     & (0.08)     & 95     \\
		Jun      & 19.59      & (0.03)     & 18.59     & (0.08)     & 95     \\
		Jul      & 16.65      & (0.03)     & 15.92     & (0.08)     & 96     \\
		Aug      & 23.11      & (0.03)     & 22.22     & (0.08)     & 96     \\
		Sep      & 24.21      & (0.03)     & 23.18     & (0.08)     & 96     \\
		Oct      & 21.49      & (0.03)     & 20.62     & (0.08)     & 96     \\
		Nov      & 19.75      & (0.03)     & 18.98     & (0.08)     & 96     \\
		Dec      & 15.40      & (0.03)     & 14.69     & (0.08)     & 95     \\ \hline
	\end{tabular}
\end{table}
\endgroup

Our benchmark is the regress-later LSM approach applied by \citet{nadarajah2018least}. This LSM version employs Monte Carlo simulation and regression to compute a VFA, which can be used to estimate a dual bound, a greedy policy, and a corresponding lower bound. Our implementation of this LSM variant is based on the same VFA specification and sample paths of the vectors of forward curves employed by the CBCD method. \looseness=-1

We utilize the same set of 500,000 independent vectors of forward curves sample paths to estimate both the LSM- and PO-based bounds.

Table \ref{Table 1} reports the dual bound estimates obtained when using both LSM and PO, as well as both their respective standard errors and percentage ratios. The standard errors of the reported estimates have equal orders of magnitude, but the LSM-based estimated dual bounds have precisions that are more than twice compared to the PO-based ones. However, in both cases the standard errors are at most 0.6\% of their respective estimates. All the estimated PO-based dual bounds are smaller than the LSM-based ones. Their ratios vary from 95\% to 96\% and equal 95\% on average. These results suggest that obtaining dual penalties using the CBCD algorithm rather than LSM is beneficial. On average, the PO-based dual bound estimates are $9\%$ smaller than the ones obtained by using zero dual penalties and the same set of sample paths. Thus, estimating good dual bounds in the considered instances is not straightforward. \looseness=-1
\begin{figure}[h!]
	\centering
	\includegraphics[width=\textwidth]{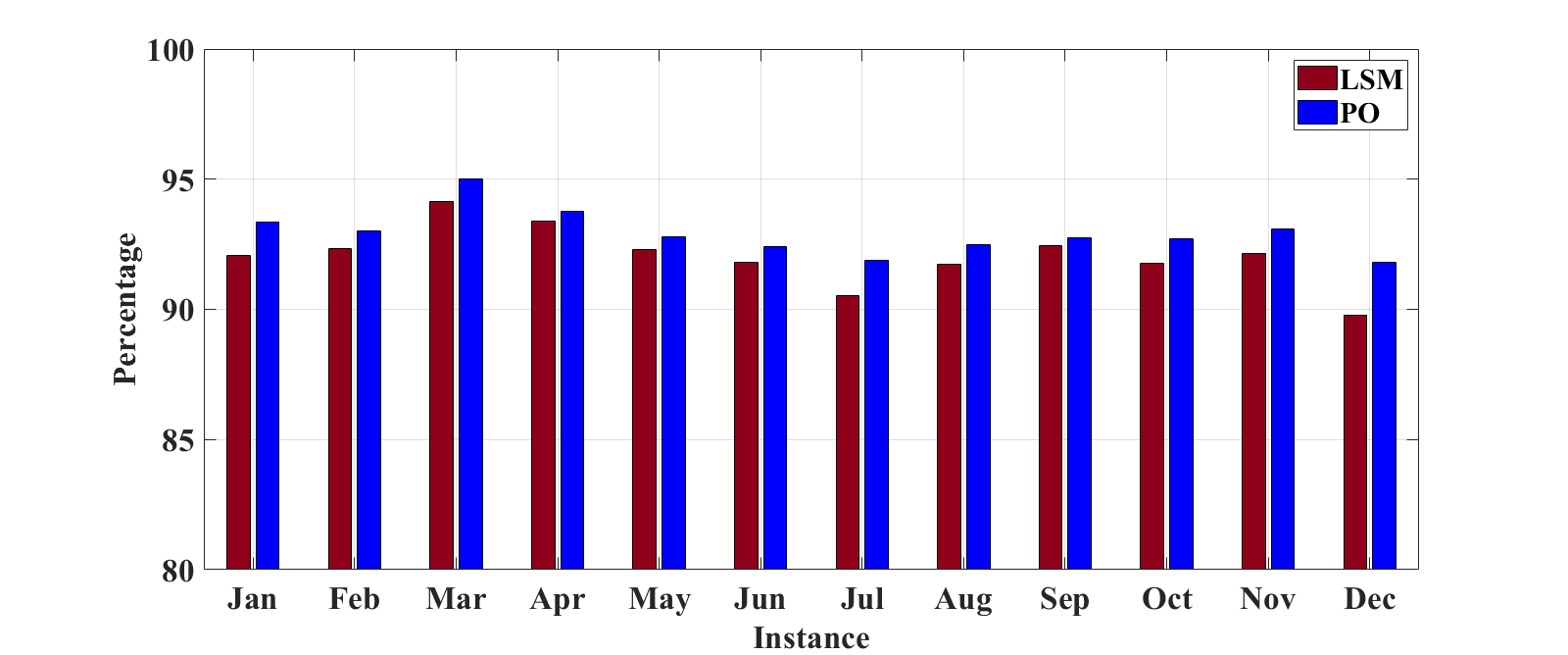}
	\caption{Comparison of the estimated LSM- and PO-based greedy policy optimality gaps as percentages of their corresponding PO-based dual bound estimates.}
	\label{Figure 4}
\end{figure}

Figure \ref{Figure 4} displays the estimated optimality gaps of the LSM- and PO-based greedy policies expressed as percentages of their respective PO-based dual bound estimates (the standard errors of their corresponding lower bound estimates are at most 0.4\% of these estimated dual bounds). The range and averages of these gaps are 6-10\% and 5-8\% and 8\% and 7\% for LSM and PO, respectively. Thus, both the PO- and LSM-based greedy policies are near optimal, but the former ones marginally outperform the latter ones. These findings corroborate the conjecture of \citet{nadarajah2018least} that the LSM-based greedy policy and dual bound are respectively near optimal and weak for the considered instances. The optimal static policy obtained by solving a deterministic version of MDP \eqref{eq:MDP} formulated based on information available in the initial stage and state has essentially zero value irrespective of the instance. Thus, using a good dynamic policy is critical in these instances.\looseness=-1	
\begin{figure}[h!]
	\centering
	\includegraphics[height=0.294\textwidth, trim={150 0 0 0} ]{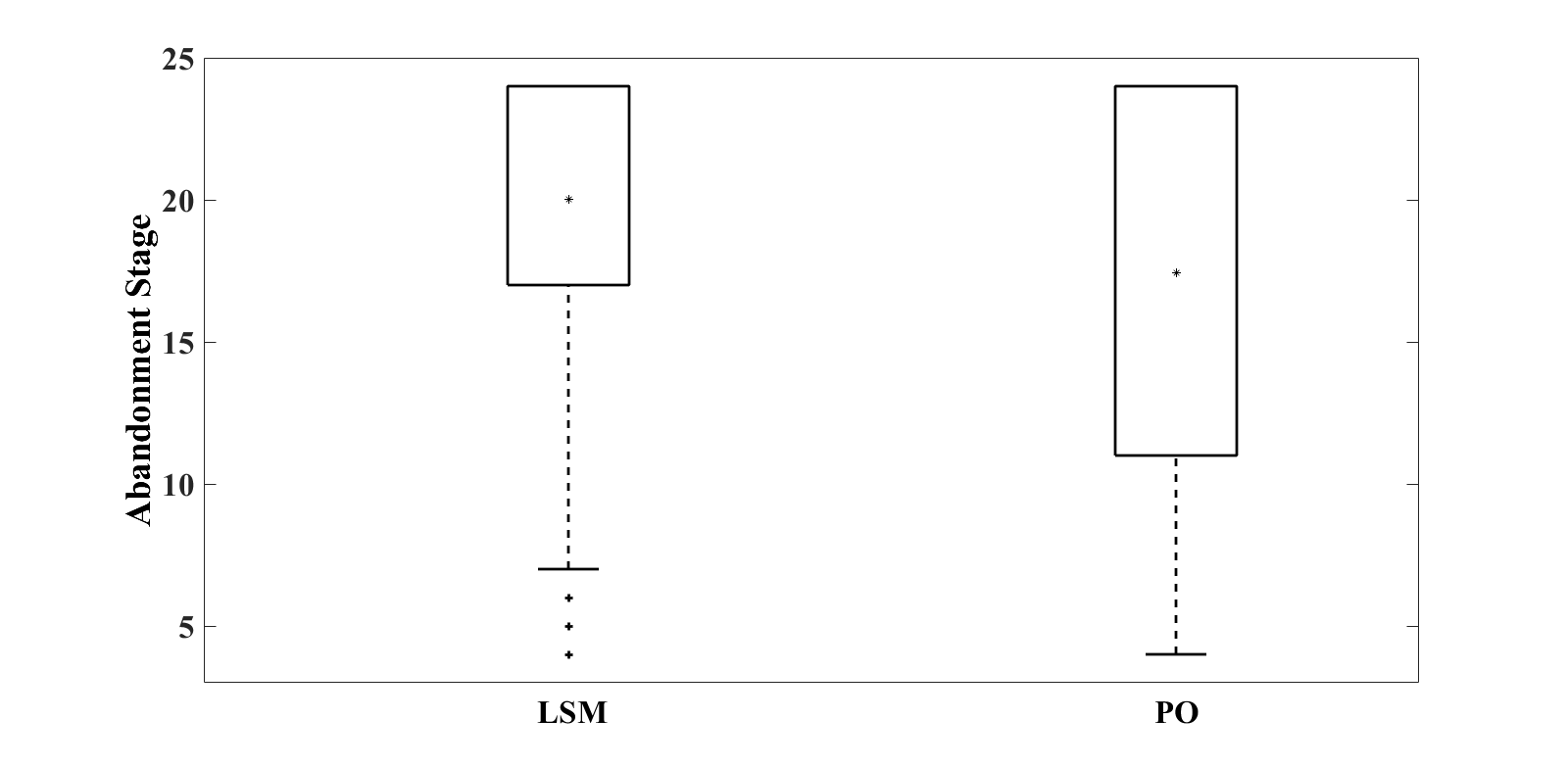}
	\caption{Frequency distributions of the stage when the LSM- and PO-based greedy policies make the abandonment decision on the January instance. The box plots show the minimum, mean, 25-th percentile, 75-th, and maximum (the 75-th percentiles and the maximum coincide).}
	\label{fig:policyAbandonStats}
\end{figure}

\begin{figure}[h!]
    \centering
\includegraphics[width=0.55\textwidth,trim={100 0 0 0}]{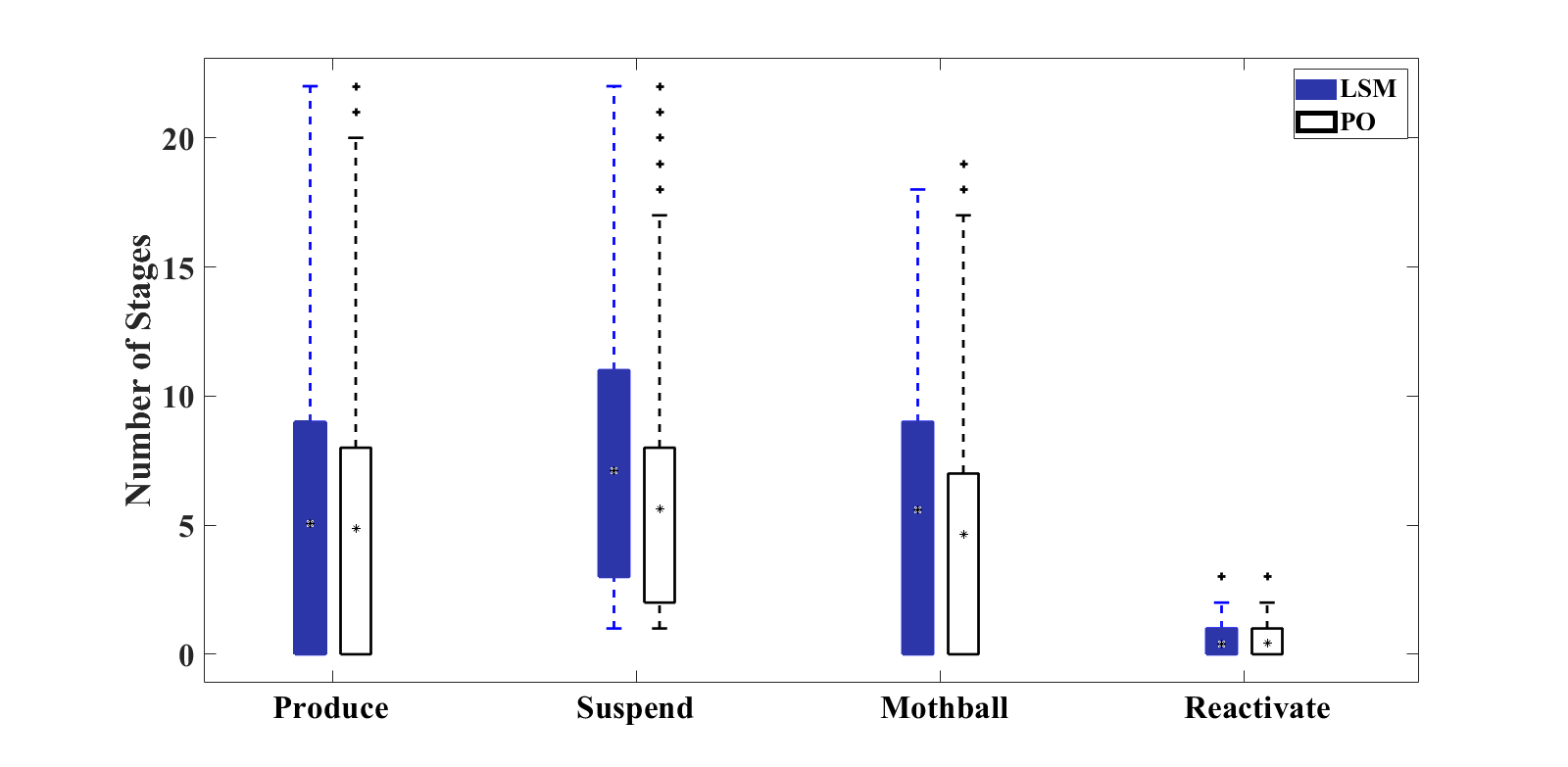}
    \caption{Frequency distributions of the non-abandonment decisions of the LSM- and PO-based greedy policies on the January instance. The box plots show the minimum, mean, 25th percentile, 75th percentile, and maximum (the minimum and the 25th percentile coincide except for the suspension action).}
    \label{fig:policyNonabandonStats}
\end{figure}
We investigate the behaviors of the LSM- and PO-based greedy policies on the representative January instance. Figure \ref{fig:policyAbandonStats} displays the frequency distributions of the stage in which these policies abandon the plant. The PO-based policy does so sooner than the LSM-based one. Figure \ref{fig:policyNonabandonStats} presents the frequency distributions of the decisions of these policies with the exception of the abandonment action. Compared to the LSM-based policy, the PO-based one produces, suspends production, and mothballs the plant fewer times, which is likely a consequence of the discussed discrepancy in their use of the abandonment decision. Thus, the LSM- and PO-based greedy policies differ in how they manage the plant's real optionality even if their values are similar. 

\begin{table}[h!]
\caption{PO and LSM average CPU times in minutes.}
\centering
\begin{tabular}{cccccccccccc}
\hline
 \multicolumn{5}{c}{PO}&\hspace{0.1in} & \multicolumn{4}{c}{LSM}\\
 \cline{1-5}\cline{7-10}
& &{Dual}&Lower&&&& {Dual}& Lower &\vspace{-0.1in}\\
 CBCD &Regression&Bound&Bound&Total&&Regression&Bound&Bound&Total\\
  \cline{1-5}\cline{7-10}
638&11&5&4&658&&5&5&4&14\\
\hline
\end{tabular}
\label{CPUTime}
\end{table}

Table \ref{CPUTime} reports the average CPU times associated with PO and LSM. We implement these methods in C++ using the GCC 4.8.5 (Red Hat 4.8.5-11) compiler and CentOS Linux 7 operating system. We use Gurobi 7.5 to solve linear programs. We apply the dlib C++ machine learning package and LAPACKE to perform PCAs and regressions, respectively. We execute our algorithms on a server with 128 GB of RAM and 12 Intel(R) Core(TM) i7-5930K processors, of which we employ at most six when running Gurobi to reduce its memory requirement. CBCD takes about ten hours to complete, which overshadows the eleven and five-and-four minutes required by the PO approach to obtain the VFAs using regression and estimate the dual-and-lower bounds. The total PO computational burden is thus roughly eleven hours. In contrast, this figure is fourteen minutes for LSM, because this technique is not based on CBCD (its regression and bound estimation efforts are comparable to the ones for PO). Despite the longer time exhibited by this algorithm, its only current alternative is to adopt a computer that has more memory than ours and attempt to solve P2LP directly using a commercial solver. \looseness=-1





\section{Conclusions}
\label{Sec:Conc}

We investigate a compound switching and timing option model of merchant energy production that gives rise to an intractable MDP. The application of current ADP techniques to realistic instances of this model provides operating policies that exhibit substantial optimality gaps. The extant literature ascribes this observation mainly to the weakness of the dual bound. We analyze this issue by applying PO to merchant energy production, extending the reach of this methodology beyond stopping models. We develop novel PCA and BCD methods to deal with the ill conditioned and large scale nature of the resulting PLP, which is out of direct reach even for a state-of-the-art commercial solver as Gurobi. Compared to LSM, despite its substantially longer run time PO leads to considerably and slightly better optimality gaps and policies, respectively, on a set of existing benchmark instances. This finding includes that both the PO- and LSM-based policies perform alternatively for the kind of models that we consider. Our research may be relevant for other commodity merchant operations situations.
\begin{APPENDIX}{Proofs}

\proof{Proof of Proposition \ref{PLPProp}.}
Let $U^{l,0}_i(x_i)$ be the value function for stage $i$ and state $x_i$ of the dynamic program \eqref{eq:samplePathDualDP2}-\eqref{eq:samplePathDualDP3} for sample path $l$ formulated with the $\beta$ vector set equal to zero.
Denote by $U^{l,0}_0(x_0)$ the corresponding value of $U^{l,\beta}_0(x_0)$.
This term is finite because it is the discounted sum of bounded rewards from the initial stage through the final one along the given sample path.
The pair $(\beta,U)$ associated with this particular choice of $\beta$ vector and the resulting $U$ vector is a feasible PLP solution with finite objective function value.
Thus, the optimal value of the PLP objective function is bounded from above.
The PLP constraint for each tuple $(l,i,x_i,\mathsf{A})$ is ${U}_{i}^{l}(x_{i})\geq r(x_i, s_i^l, \mathsf{A})$.
The right hand side of this inequality evaluates to 0 or $S$ when $x_i$ equals $\mathsf{A}$ or it belongs to $\mathcal{X}^{\prime}$.
It follows that the optimal PLP objective function value is bounded from below.
Suppose that all optimal PLP solutions have at least one infinite element of their corresponding $\beta$ vector.
Pick an arbitrary optimal PLP solution $(\beta^\ast,U^\ast)$.
Let $(0,U^{(0)})$ be the PLP solution used to establish that the optimal PLP objective function value is bounded from above.
In particular, it is a basic feasible solution for PLP.
Consider a sequence of such solutions that starts from $(0,U^{(0)})$, ends at $(\beta^\ast,U^\ast)$, and includes as additional elements, if any, points that belong to the boundary of the PLP feasible set.
The penultimate item of this sequence is a vertex that has both a ray that connects it to $(\beta^\ast,U^\ast)$ and zero objective function gradient, because the $\beta^\ast$ vector has at least one infinite element, the solution $(0,U^{(0)})$ is finite, and the optimal PLP objective value is bounded from both below and above.
That is, it is a finite optimal PLP solution, which contradicts the assumption that all optimal PLP solutions have at least one infinite element.
Thus, PLP has at least one finite optimal solution.
\Halmos
\endproof

\proof{Proof of Proposition \ref{prop:equivPCAPLP}.}
Pick a feasible PLP solution $(\beta,U)$.
Define $\beta^\prime$ as the vector with $(i,x_i)$ component $\beta^\prime_{i,x_i}$ equal to $W^{-1}_{i,x_i} \beta_{i,x_i}$.
Evaluating the left hand sides of the PLP constraints and the P2LP ones at $(\beta,U)$ and $(\beta^\prime,U)$, respectively, yields the same values.
An analogous result holds for the feasible P2LP solution $(\beta,U)$ and the PLP one $(\beta^\prime,U)$ for which the $(i,x_i)$ part $\beta^\prime_{i,x_1}$ of the vector $\beta^\prime$ is $W_{i,x_i} \beta_{i,x_i}$.
That is, there is a one to one mapping between the respective sets of PLP and P2LP feasible solutions. 
PLP and P2LP have the same objective function.
Thus, their optimal solution sets coincide.
\Halmos
\endproof

\proof{Proof of Proposition \ref{PBCDConverges}.} 
The Bolzano-Weierstrass theorem and Proposition 2 of \cite{grippo2000convergence} imply that the sequence of solutions generated by the idealized CBCD algorithm converges.
Denote by $\overline{\beta}$ the part of the resulting solution that corresponds to the $\beta$ vector of decision variables.
Suppose that for each sample path $l\in\mathcal{L}$ of the vector of forward curves the linear program \eqref{DPLPObj}-\eqref{DPLPConst3} formulated with $\beta$ equal to $\overline{\beta}$ has a non-degenerate optimal solution vector $U^{l,\ast}(\overline{\beta})$.
This assumption implies that its dual model has a unique optimal solution, which we denote as $\mu^{l,\ast}(\overline{\beta})$.
Define $U^\ast(\overline{\beta}):=(U^{l,\ast}(\overline{\beta}),l\in\mathcal{L})$ and $\mu^\ast(\overline{\beta}):=(\mu^{l,\ast}(\overline{\beta}),l\in\mathcal{L})$.
The pair $(\overline{\beta},U^\ast(\overline{\beta}))$ is a feasible solution for P2LP.
Consider the idealized P2LP dual model.
It has (i) vectors of decision variables $\mu$ and $\theta$ that are associated with \eqref{form:Const} and \eqref{betaBounds} and (ii) two sets of constraints that are related to the $\beta$ and $U$ vectors of the variables of the idealized PLP.
The $\beta$-related constraints and the $\theta$ variables can be subdivided according to the elements of the given partition $\mathcal{P}$.
At each iteration, the idealized CBCD method solves the linear program \eqref{form:ObjBCDLP}-\eqref{betaBounds}, which we denote as LP$^p$, for each set $\mathcal{P}_p\in\mathcal{P}$.
The dual of this model features the $\mu$ and $\theta(\mathcal{P}_p):=(\theta_{i,x_i,b}, (i,x_i,b) \in \mathcal{P}_p\times\mathcal{B}_i)$ variable vectors and the $U$- and $\beta(\mathcal{P}_p)$-related constraints.
The pair $(\overline{\beta}(\mathcal{P}_p),U^\ast(\overline{\beta}))$ is an optimal LP$^p$ solution.
Consider the complementary slackness conditions for this model and its dual expressed with respect to it.
Solving the ones associated with $U^\ast(\overline{\beta})$ amounts to finding a solution to the system of $U$-related constraints that define $\mu^\ast(\overline{\beta})$, which we know uniquely solves it.
Moreover, because we assume that $\overline{\beta}(\mathcal{P}_p)$ strictly satisfies the inequalities \eqref{betaBounds}, the elements of the corresponding optimal vector $\theta^\ast(\mathcal{P}_p; \overline{\beta}(\mathcal{P}_p))$ equal zero.
Thus, $(\overline{\beta},U^\ast(\overline{\beta}))$ and $(\mu^\ast(\overline{\beta}),\theta^\ast(\mathcal{P}_p; \overline{\beta}(\mathcal{P}_p)))\equiv (\mu^\ast(\overline{\beta}),0)$ comply with both the complementary slackness equations for the idealized P2LP and its dual and the ones for P2LP and its dual.
Moreover, the pair $(\mu^\ast(\overline{\beta}),0)$ fulfills the set of $\beta(\mathcal{P}_p)$-related constraints of the idealized P2LP dual for each $p\in\{1,\ldots,P\}$, because it is the only solution for the dual of LP$^p$ that satisfies complementary slackness with respect to its primal optimal solution $(\overline{\beta}(\mathcal{P}_p),U^\ast(\overline{\beta}))$.
It follows that $\mu^\ast(\overline{\beta})$ is feasible for the dual of P2LP.
Consequently, $(\overline{\beta},U^\ast(\overline{\beta}))$ and $\mu^\ast(\overline{\beta})$ optimally solve P2LP and its dual, respectively.
\Halmos
\endproof

 \end{APPENDIX}

\ACKNOWLEDGMENT{%
This research is in part supported by a grant from the Scott Institute for Energy Innovation at Carnegie Mellon University and by NSF grant CMMI 1761742. The second author acknowledges
support from the UIC College of Business Dean’s Research Grant.
The third author is a Faculty Affiliate of the Scott Institute for Energy Innovation at Carnegie Mellon University.
}



\bibliographystyle{informs2014}
\bibliography{references} 
\end{document}